\documentclass{article}    
\usepackage{graphicx}
\usepackage{color}
\usepackage{multirow}
\usepackage{graphics} 
\usepackage{epsfig} 
\usepackage{amsmath} 
\usepackage{amssymb}  
\usepackage{amsfonts}

\usepackage{overpic}

\DeclareMathOperator{\aew}{\forall}

\newtheorem{assum}{Assumption}
\newtheorem{rem}{Remark}
\newtheorem{defn}{Definition}
\newtheorem{thm}{Theorem}
\newtheorem{lem}{Lemma}
\newtheorem{cor}{Corollary}
\newtheorem{prop}{Proposition}

\title{Backward-Forward Reachable Set Splitting for State-Constrained Differential Games}
\author{Xuhui Feng, Mario E.~Villanueva, Boris Houska}
\date{}

\begin{document}

\maketitle

\begin{abstract}
This paper is about a set-based computing method for solving a general class of 
two-player zero-sum Stackelberg differential games. We assume that the game is 
modeled by a set of coupled nonlinear differential equations, which can be influenced
by the control inputs of the players. Here, each of the players has to satisfy 
their respective state and control constraints or loses the game. The main 
contribution is a backward-forward reachable set splitting scheme, which can be 
used to derive numerically tractable conservative approximations of such two 
player games. In detail, we introduce a novel class of differential inequalities
that can be used to find convex outer approximations of these backward and forward
reachable sets. This approach is worked out in detail for ellipsoidal set 
parameterizations. Our numerical examples illustrate not only the effectiveness 
of the approach, but also the subtle differences between standard robust optimal
control problems and more general constrained two-player zero-sum Stackelberg 
differential games.\\[0.2cm]
\end{abstract}

\section{Introduction}
The origins of differential games and game-theoretic optimal control go back
to~\cite{Isaacs1999}. A historical overview of the early developments, roughly
ranging from 1950--1970, can be found in~\cite{Breitner2005}. There it becomes
clear that, since their inception, the theories of differential games and optimal
control have been deeply intertwined~\cite{Pesch2009}. The mathematical foundation
of modern differential game theory, was established between 1970-1990. Precise
definitions and a mature mathematical framework for this theory can, for example,
be found in the books by~\cite{Friedman1971} as well as~\cite{Krasovskii1988}.

The question of how to define an appropriate mathematical model for a game has,
in general, no unique answer. This is due to the fact that the construction of
such a model may depend on many different aspects. One of these aspects is the 
information each player has about other players, including their goals, their 
ability (and willingness) to communicate, and their willingness to agree on 
actions and rules~\cite{Basar1986,Basar1999,Nash1951}. The focus of this paper 
is on a specific class of games, namely, two-player zero-sum Stackelberg 
games~\cite{Bensoussan2018,Stackelberg1952}, also known as 
``worst-case games''. Here, Player~1 makes a decision and announces it to the 
other player, who has a conflicting objective. Stackelberg games arise in economics, marketing and policy making~\cite{Dockner2000}, but
their applicability also extends to control systems~\cite{Basar1979,Papavassilopoulos1979}.

In a context of mathematical programming applied to static zero-sum Stackelberg
games, one can distinguish between semi-infinite programming (SIP)~\cite{Hettich1993}
and generalized semi-infinite programming (GSIP) problems~\cite{Jongen1998}. 
Both classes of problems can be used to model zero-sum Stackelberg games, but in SIP
the feasible set of the second player is assumed to be constant. This is in contrast to GSIP, where the decision of the first player affects the feasible set of the optimization problem of the second player. Notice that there exists a vast body of literature on numerical methods for both SIP and GSIP~\cite{Stein2012,Diehl2013}.

The main difference between static and differential games is that, in the latter
case, the state of the game is modeled by a differential equation~\cite{Friedman1971},
\begin{equation*}
\forall t \in [0,T], \quad \dot x(t) = f(x(t),u_1(t),u_2(t)) \; .
\end{equation*}
Here, the decision variables of the players are functions of time: Player~1
chooses the input function $u_1$, while Player~2 chooses $u_2$. A Stackelberg 
differential game can either be played in open-loop or closed-loop mode. In the 
former, Player~1 chooses $u_1$ first and announces his decision to Player~2. In 
contrast, in closed-loop mode both players make their decisions simultaneously
and continuously in time, for example, based on the current state $x(t)$ of the system.

Worst-case robust optimal control problems are a special class of 
zero-sum differential games, where the control $u_2$ is bounded by a given set. 
For example, $H_\infty$ control problems can be interpreted as differential games
against nature~\cite{Bacsar2008}. In the mathematical control literature robust
optimal control problems are frequently analyzed by means of Hamilton-Jacobi 
equations~\cite{Bardi1991,Bressan2006}. Modern numerical methods for robust
optimal control, both in open- and closed-loop, are often based on
set-theoretic considerations~\cite{Blanchini2008,Houska2018,Langson2004}.

In order to understand the contributions of this paper, it is important to be
aware of one fact: standard robust optimal control problems are differential games
which do not enforce state constraints on the adverse player. This is in contrast
to differential games played against a rational opponent, who 
might develop sophisticated strategies but agrees to act subject to the rules of
the game. Formally, such rules can be modeled via state constraints
$\mathbb X_1(t) \subseteq \mathbb R^{n_x}$ and $\mathbb X_2(t) \subseteq \mathbb R^{n_x}$. 
These sets are such that Player~1 loses if the constraint 
\begin{equation}
\notag
\forall t \in [0,T], \qquad x(t) \in \mathbb X_1(t)
\end{equation}
is violated, while Player~2 loses if the constraint
\begin{equation}
\notag
\forall t \in [0,T], \qquad x(t) \in \mathbb X_2(t)
\end{equation}
is violated. For example, in a football game an offside can be modeled by a
suitable state constraint set $\mathbb X_2(t)$ that is related to the relative 
positions (differential states) of the players and the ball. Rational players could,
in principle, exploit their explicit knowledge of the state constraints of their
opponents. In football, for example, one team could let the adverse team run into
an offside trap. To draw an analogy to static games, one could state that 
zero-sum Stackelberg differential games with state constraints are to robust optimal
control problems what GSIP problems are to SIP problems. In the context of general
differential games, the actions of the first player affect the set of feasible 
actions of the second player. This is because the solution of the differential 
equation depends on both $u_1$ and $u_2$.

Compared to the vast amount of literature on numerical methods for standard robust
optimal control, the number of articles on zero-sum differential games with state
constraints on both players is rather limited. An interesting historical example 
for a differential game with state constraints is the famous ``man and lion'' problem. In this game, both players have equal maximum speed, and are both constrained to stay
in a circular arena. The rather surprising fact that the man can survive infinitely
long without being caught by the lion has been proven in 1952 by Besicovitch;~\cite[see]{Littlewood1986}.
Most contributions in the area of state-constrained zero-sum differential games use, in one way or another, concepts from viability 
theory~\cite{Aubin2009}. For example, \cite{Cardaliaguet2001} used viability 
kernel techniques to construct numerical methods for differential games with 
separable dynamics and state constraints.  An overview of recent 
advances in zero-sum differential games with state constraints can also be found
in~\cite{Cardaliaguet2007}.

The main contribution of this paper is a set-based computing framework for
analyzing and constructing approximate, yet conservative solutions of zero-sum
Stackelberg differential games. The proposed framework is able to deal with 
coupled dynamics as well as state constraints for both players. The set-based 
problem formulation is outlined in Section~\ref{sec::setting}. In contrast
to~\cite{Cardaliaguet2001}, we do not assume that the dynamics of the game are 
separable. Section~\ref{sect::bfRsS} introduces a generic backward-forward reachable
set splitting result for such non-separable two player differential games, which
is presented in Theorem~\ref{thm::splitreach}. This result is then used to 
construct convex outer approximations of the constrained reachable set of the 
second player via a system of generalized differential inequalities, as summarized
in Theorem~\ref{thm::outerreach}. Theorem~\ref{thm::ellreach} specializes
this result for ellipsoidal set parameterizations. In Section~\ref{sec::case}, the latter construction is leveraged in order to construct a standard optimal control problem with boundary constraints, whose solution conservatively approximates the
solution of the original set-based problem. This problem can be solved using 
existing, gradient-based, optimal control algorithms. The developments of 
this paper are demonstrated through a numerical example implemented using the 
optimal control software {\tt ACADO Toolkit}~\cite{Houska2011}, in Section~\ref{sec::numerics}. Section~\ref{sec::conclusion} concludes the paper.

\paragraph*{Notation} The set of $n$-dimensional $L_1$-integrable functions is 
denoted by $\mathbb L_1^n$ while $\mathbb{W}^{n}_{1,1}$ denotes the associated
Sobolev space of weakly differentiable functions with $L_1$ integrable derivatives.
The set of compact and convex compact subsets of $\mathbb R^n$ are denoted by 
$\mathbb K^n$ and $\mathbb K^n_C$, respectively. For a set $Z\subseteq\mathbb{R}^{n}$,
$\mathcal{P}(Z)\subseteq\mathbb{R}^{n}$ denotes its power set, which includes the
empty set, denoted by $\varnothing$. Moreover, $\operatorname{int}(Z)$ denotes the
interior of a set $Z\subseteq\mathbb{R}^{n}$ and $\operatorname{cl}(Z)$ its closure
in $\mathbb R^n$. The support function of a set $Z$ is defined
as
\begin{equation*}
\forall c\in\mathbb{R}^{n},\quad 
V[Z](c) = \sup_{z\in Z} \ c^\intercal z \;.
\end{equation*}
Additionally, we define $V[\varnothing](c) = -\infty$.
The sets of positive semidefinite and positive definite $n$-dimensional
matrices is denoted by $\mathbb{S}^{n}_{+}$ and $\mathbb{S}^{n}_{++}$. 
An ellipsoid with center $q \in \mathbb R^n$ and shape matrix $Q \in \mathbb S^n_+$
is given by 
\begin{equation*}
\mathcal E(q, Q) = \left\{ 
q + Q^{\frac{1}{2}} v \ \middle | \ \exists v\in\mathbb{R}^{n}: \
v^\intercal v \leq 1
\right\} \; ,
\end{equation*}
where $Q^{\frac{1}{2}}$ can be any square root of $Q$, as the unit
disc in $\mathbb R^{n}$ remains invariant under orthogonal transformations.

\section{Open-loop Stackelberg differential games}
\label{sec::setting}

This paper is about two-player differential games whose state, 
$x \in \mathbb{W}^{n_x}_{1,1}$, satisfies a differential equation of 
the form
\begin{equation}
\label{eq::odes}
\begin{alignedat}{3}
&\aew t\in[0,T],\quad
&&\dot{x}(t) &&= f( x(t), u_1(t), u_2(t) ),  \\
&\text{with} &&x(0) &&= x_{0} \;.
\end{alignedat}
\end{equation}
Here, Player 1 chooses the control input $u_1: [0,T] \to \mathbb U_1$ while 
Player 2 chooses the control input $u_2: [0,T] \to \mathbb U_2$.
\begin{assum}
\label{ass::rhs1}
The right-hand side function, $f$, is jointly continuous in $x,u_1,u_2$ and 
locally Lipschitz continuous in $x$.  
\end{assum}
\begin{assum}
\label{ass::admissible}
The control constraint sets $\mathbb U_1, \mathbb U_2 \subseteq \mathbb R^{n_u}$ are non-empty, convex, and compact.
\end{assum}
For simplicity of presentation, it is also assumed that the initial state 
$x_0 \in \mathbb R^{n_x}$ is a given constant. The developments in this paper can
easily be generalized for the case that the initial value is chosen
by one of the players or to the case where the inputs have different
dimensions.

\begin{rem}
\label{rem::1}
Two-player differential games can sometimes be represented using a differential 
equation system of the form
\begin{equation}
\label{eq::odes2}
\begin{alignedat}{3}
&\aew t\in[0,T],\quad
&&\dot{x}_1(t) &&= f_1( x_1(t), u_1(t) ),  \\
&\aew t\in[0,T],\quad
&&\dot{x}_2(t) &&= f_2( x_2(t), u_2(t) ),  \\
&\text{with} &&x(0) &&= x_{0} \; .
\end{alignedat}
\end{equation}
These systems are called separable, since each function $f_i$, $i \in \{ 1,2\}$, 
depends only on the state $x_i$ and control $u_i$. As an example, consider a 
simple two-car race: the states of each car (position, orientation, and velocities)
are only functions of its own controls (acceleration and steering). Likewise, there 
are systems that cannot be formulated as~\eqref{eq::odes2} easily. 
Consider for example two children on a seesaw: each child may shift its own 
weight independently (control input), while the state of the system (position and
velocity of the seesaw) is simultaneously affected by both control inputs. We keep
the formulation general, as every separable system can be written as~\eqref{eq::odes}
by introducing the stacked state $x^\intercal = (x_1^\intercal ,x_2^\intercal)$.
\end{rem}

\subsection{State constraints}
State constraints can be used to define the rules of a game. Here, we consider 
two given---and potentially time varying---state constraint sets, $\mathbb X_1(t),\mathbb X_2(t)
\subseteq \mathbb R^{n_x}$. Player~1 loses if the constraint 
\begin{equation}
\label{eq::player1}
\forall t \in [0,T], \qquad x(t) \in \mathbb X_1(t)
\end{equation}
is violated. Likewise, Player~2 loses if the constraint
\begin{equation}
\label{eq::player2}
\forall t \in [0,T], \qquad x(t) \in \mathbb X_2(t)
\end{equation}
is violated. Depending on the particular definitions of $\mathbb X_1$ and 
$\mathbb X_2$, there may be situations in which both players lose, one of the 
players loses, or no one loses.

\begin{rem}
There are games, where the state constraint sets $\mathbb X_1$ and $\mathbb X_2$
coincide. Consider again the two-car race from Remark~\ref{rem::1}: one may be interested in enforcing a collision avoidance constraint. This gives rise to a coupled state constraint involving the positions of both cars. If a collision occurs, both cars
are out of the race. Thus, if there are no further constraints, we have 
$\mathbb X_1 = \mathbb X_2$. However, as soon as we introduce the additional rule
that Player~1 loses if the first car leaves the road while Player~2 loses if the
second car is not staying on track, we have $\mathbb X_1 \neq \mathbb X_2$ (in 
this example, Player~1 does not necessary lose if the second player's car is not
staying on the road and vice-versa).
\end{rem}

\subsection{Feasibility}

We use the symbol $\mathbb{X}[u_1]$ to denote the set of feasible state trajectories
that the second player can realize,
\begin{equation}
\label{eq::reachset}
\mathbb{X}[u_1] = \left\{
x \in \mathbb W^{n_x}_{1,1} \ \middle |
\begin{array}{l}
\exists u_2 \in \mathbb L^{n_u}_{1}: \ \aew \tau \in [0,\,T], \\
\dot x(\tau) = f(x(\tau),u_1(\tau),u_2(\tau)), \\
x(\tau) \in \mathbb{X}_2(\tau), \ u_2(\tau) \in \mathbb{U}_2, \\
x(0)  = x_0
\end{array}
\right\}.
\end{equation}
Thus, the set-valued function $X[u_1]:\mathbb{R}\to
\mathcal{P}\left(\mathbb{R}^{n_x}\right)$, given by
\begin{equation*}
\forall t \in [0,T], \qquad X[u_1](t) = \left\{ \ x(t) \in \mathbb R^{n_x} \ 
\middle | \ x \in \mathbb{X}[u_1] \ \right\}\;,
\end{equation*}
denotes the reachable set in the state space.

\begin{defn}
\label{def::lowerF}
A control input $u_{1}\in\mathbb{L}^{n_u}_{1}$ with $u_1: [0,T] \to \mathbb U_1$
is called lower-level feasible if
\begin{equation*}
\mathbb{X}[u_1] \neq \varnothing \; .
\end{equation*}
Otherwise, $u_1$ is called lower-level infeasible.
\end{defn}

Notice that if Player~1 chooses a control input $u_1$ that is lower-level 
infeasible, Player~2 is forced to violate the rules of the game. Because this paper
focuses on games in which such behavior of Player~1 is unwanted, we introduce the
following definition of upper level feasibility.

\begin{defn}
\label{def::upperF}
A control input $u_1\in\mathbb{L}^{n_u}_{1}$ with $u_1: [0,T] \to \mathbb U_1$ is
called upper-level feasible if it is lower-level feasible and 
\begin{equation*}
\forall t\in[0,T], \quad X[u_1](t) \subseteq \mathbb X_1(t) \; .
\end{equation*}
Otherwise, $u_1$ is called upper-level infeasible.
\end{defn}

\begin{rem}
Enforcing the constraint $\mathbb{X}[u_1]\neq \varnothing$ is equivalent to introducing the 
rule that Player~1 loses the game if $u_1$ is lower-level infeasible. At this 
point, one should be clear in mind that requiring lower-level feasibility does 
not imply that Player~1 is not allowed to ``win'' the game. The conditions for 
lower---and upper---level feasibility merely define under which conditions 
Player~1 loses the game.
\end{rem}

\begin{rem}
Robust optimal control~\cite{Houska2012} considers games with $\mathbb X_2 = \mathbb R^{n_x}$
with the disturbances being the input of the adverse player. In such
case, all inputs $u_1$ are lower-level feasible, since the second player has no state 
constraints that could possibly become infeasible.
\end{rem}

\subsection{Constrained open-loop zero-sum differential games}

The goal of this paper is to analyze and approximately solve constrained 
open-loop zero-sum differential games with Stackelberg information structure. 
Thus, we assume that Player~1 chooses a strategy and announces it 
to Player~2. An optimal open loop strategy for Player~1 is any solution of 
\begin{equation}
\label{eq::exactgame}
\begin{aligned}
&\inf_{u_1} \  \mathcal M( X[u_1](T) ) \\
&\text{s.t.} \ \left\lbrace
\begin{alignedat}{3}
&\mathbb{X}[u_1] &&\neq \varnothing &&\\
&X[u_1](t) &&\subseteq \mathbb{X}_{1}(t) \quad &&\text{for all} \; \; t\in[0,T]\\
& u_1(t) &&\in \mathbb{U}_1  &&\text{for all} \; \; t\in[0,T] \; .
\end{alignedat} 
\right.
\end{aligned}
\end{equation}
We assume that a Mayer term $m: \mathbb R^{n_x} \to \mathbb R$ is given and
that
\begin{equation*}
\mathcal M(X) = \sup_{\xi \in X} \ m( \xi)
\end{equation*}
denotes the supremum (worst-case value) of $m$ on a given set $X$. Notice 
that~\eqref{eq::exactgame} is feasible if and only if there exists a control input
$u_1$ that is upper-level feasible. Additionally, we recall that if $u_1$ is a 
feasible point of~\eqref{eq::exactgame}, then $u_1$ \mbox{is---by} construction of the 
rules of our \mbox{game---also} lower-level feasible. 

\begin{rem}
If there exists an equilibrium solution $(u_1^\star,u_2^\star)$ (min-max point) of the open-loop Stackelberg zero-sum
differential game, it can be computed by first finding a minimizer $u_1^\star$ of~\eqref{eq::exactgame}. In this case, $u_2^\star$ must be a maximizer of
\begin{equation*}
\max_{x,u_2} \ m(x(T)) \quad \mathrm{s.t.} \
\begin{cases}
\forall t\in [0,T]: \\
\dot x(t) = f(x(t),u^\star_{1}(t),u_2(t))\;,  \\
x(t) \in \mathbb{X} \;, \ u_2(t) \in\mathbb{U}_{2} \\
x(0) = x_0 \; .
\end{cases} 
\end{equation*}
\end{rem}

Throughout this paper it is assumed that the decision of whether one of the 
players has won the game is made a posteriori, after both players have 
played their strategies. This final decision is based on the objective value 
$m(x(T))$, where $x(T)$ denotes the state of the system at time $T$. As we can 
always add constant offsets to $m$, we say that Player~1 wins the game if $u_1$ 
is upper-level feasible and $m(x(T)) \leq 0$. Similarly, Player~2 wins if 
$x(t) \in \mathbb X_2(t)$ for all $t \in [0,T]$ and $m(x(T)) > 0$. Notice that 
this definition is consistent in the sense that at most one player can win the 
game and it is impossible that one of the players wins and loses a 
game simultaneously.

\begin{rem}
One can also consider closed-loop games whose dynamics are given by
\begin{equation*}
\label{eq::odesClosedLoop}
\begin{alignedat}{3}
&\forall t\in[0,T],\quad
&&\dot{x}(t) &&= F( x(t), \mu_1( t, x(t)), \mu_2(t, x(t) ) \\
&\text{with} &&x(0) &&= x_{0} \;.
\end{alignedat}
\end{equation*}
Here, Player~1 chooses the feedback law $\mu_1: \mathbb R \times \mathbb R^{n_x} \to \mathbb U_1$
while Player~2 chooses the feedback law $\mu_2: \mathbb R \times \mathbb R^{n_x} \to \mathbb U_2$.
Closed-loop games in full generality are significantly harder to solve and analyze
than open-loop games, even without the presence of state-constraints. 
Nevertheless, by restricting the search to parametric feedback laws, for example, affine 
feedback laws of the form
\begin{equation*}
\mu_1(x) = K_1(t) x + k_1(t) \quad \text{and} \quad \mu_2(x) = K_2(t) x + k_2(t) \; ,
\end{equation*}
one can formulate closed-loop games in the form of~\eqref{eq::exactgame} by regarding
the control law coefficients $u_1 = ( \text{vec}(K_1)^\intercal, k_1^\intercal )^\intercal$ and $u_2 = ( \text{vec}(K_2)^\intercal, k_2^\intercal )^\intercal$ as the inputs of the first and second player, respectively.
\end{rem}

\section{Backward-forward reachable set splitting}
\label{sect::bfRsS}
The goal of this section is to analyze the reachable sets $X[u_1](t)$. 
Notice that, due to the presence of state constraints for Player 2, 
it is non-trivial to ensure lower-level feasibility of an input $u_1$ (see Remark~4).
In particular, the state constraints $\mathbb{X}_2(\cdot)$ induce a coupling 
in time; that is, knowing the reachable state at time $t$, requires knowledge of
the state trajectories on the whole time horizon.

In order to remove the coupling in time, we introduce the set-valued 
function $X_{\rm B}[u_1]:\mathbb{R}
\to\mathcal{P}\left(\mathbb{R}^{n_x}\right)$ given by
\begin{equation}
\label{eq::breachset}
X_{\mathrm{B}}[u_1](t) = \left\{ \xi \in \mathbb R^{n_x}
\middle | \
\begin{aligned}
&\exists x \in \mathbb W^{n_x}_{1,1},\ \exists u_2 \in \mathbb L^{n_u}_{1}: \\
&\aew \tau \in [t,T], \\
&\dot x(\tau) = f(x(\tau),u_1(\tau),u_2(\tau))\\
&x(\tau) \in \mathbb X_2(t), \ u_2(\tau) \in \mathbb U_2, \\
&x(t)  = \xi
\end{aligned}
\right\} .
\end{equation}
The set $X_{\mathrm{B}}[u_1](t)$ can be interpreted as the set of all states $x(t)$
of the game at time $t$ for which the second player is able to satisfy the rules
of the game on the remaining time interval $[t,T]$. Next, we introduce the set-valued
function $X_{\rm F}[u_1]:\mathbb{R}\to
\mathcal{P}\left(\mathbb{R}^{n_x}\right)$ given by
\begin{equation}
\label{eq::freachset}
\begin{aligned}
&X_{\mathrm F}[u_1](t) = \left\{ \xi \in \mathbb R^{n_x}
\middle | \
\begin{aligned}
&\exists x \in \mathbb W^{n_x}_{1,1},\ \exists u_2 \in \mathbb L^{n_u}: \\
&\aew \tau \in [0,t], \\
&\dot x(\tau) = f(x(\tau),u_1(\tau),u_2(\tau))\\
&x(\tau) \in X_{\rm B}[u_{1}](\tau),\ u_2(\tau) \in \mathbb U_2, \\
&x(0)  = x_0 \ , \ x(t) = \xi
\end{aligned}
\right\}.
\end{aligned}
\end{equation}
The set-valued function $X_{\rm F}[u_1](t)$ is called the constrained 
forward reachable set of the system at time $t$. The following
theorem establishes the fact that $X$ and $X_{\rm F}$ coincide.

\begin{thm}
\label{thm::splitreach}
The equation $X[u_1] = X_{\rm F}[u_1]$ holds for all input functions $u_1 \in \mathbb L_1^{n_u}$.
\end{thm}

\textbf{Proof.}
Let $t \in [0,T]$ and $u_1 \in \mathbb L^{n_u}_{1}$ be given.
The goal of the first part of this proof is to establish the inclusion
\begin{equation}
\label{eq::inclusion1}
X_{\rm F}[u_1](t) \subseteq X[u_1](t) \; .
\end{equation}
Let the functions $x_\mathrm{F} \in \mathbb W^{n_x}_{1,1}$ and $u_{2,\mathrm{F}}: [0,t] \to \mathbb U_2$ be such that
\begin{align}
\label{eq::auxodeF}
\dot x_\mathrm{F}(\tau) &= f(x_\mathrm{F}(\tau),u_1(\tau),u_{2,\mathrm{F}}(\tau))\\
\label{eq::aux11}
x_\mathrm{F}(\tau) &\in X_{\mathrm{B}}[u_{1}](\tau),\\
\label{eq::initF}
x_\mathrm{F}(0)  &= x_0
\end{align}
for all $\tau \in [0,t]$. Now, the definition of $X_{\rm B}[u_{1}]$ 
and~\eqref{eq::aux11} imply that there exists a $x_\mathrm{B} \in \mathbb W_{1,1}^{n_x}$ and $u_{2,\mathrm{B}}: [t,T] \to \mathbb U_2$ such that
\begin{align}
\label{eq::auxodeB}
\dot x_\mathrm{B}(\tau) &= f(x_\mathrm{B}(\tau),u_1(\tau),u_{2,\mathrm{B}}(\tau))\\
\label{eq::inclB}
x_\mathrm{B}(\tau) &\in \mathbb X_2(\tau),\\
\label{eq::initB}
x_\mathrm{B}(t)  &= x_\mathrm{F}(t)
\end{align}
for all $\tau \in [t,T]$. Thus, we can construct the functions
\begin{align*}
x(\tau) &= 
\begin{cases}
x_\mathrm{F}(\tau) &\text{if} \ \ 0 \leq \tau \leq t \\
x_\mathrm{B}(\tau) &\text{if} \ \ t \leq \tau \leq T
\end{cases}
\intertext{and}
u_2(\tau) &= 
\begin{cases}
u_{2,\mathrm{F}}(\tau) &\text{if} \ \ 0 \leq \tau \leq t \\
u_{2,\mathrm{B}}(\tau) &\text{if} \ \ t \leq \tau \leq T
\end{cases}
\end{align*}
with $u(\tau)\in\mathbb{U}_2$ for all $\tau\in[0,T]$. These functions satisfy
\begin{alignat}{2}
\label{eq::odeT}
\dot x(\tau) 
&\overset{\eqref{eq::auxodeF},\eqref{eq::auxodeB},\eqref{eq::initB}}{=}
&& f(x(\tau),u_1(\tau),u_{2}(\tau))\\
\label{eq::inclT}
x(\tau) 
&\overset{\hphantom{5)}\eqref{eq::aux11},\eqref{eq::inclB}\hphantom{,(1}}{\in} 
&& \mathbb X_2(\tau),\\
\label{eq::initT}
x(0)
&\overset{\hphantom{(15)}\eqref{eq::initF}\hphantom{,(15)}}{=} 
&& x_0
\end{alignat}
Inclusion~\eqref{eq::inclT} follows from $\eqref{eq::aux11}$ and $\eqref{eq::inclB}$ 
since, by construction, $X_{\rm B}$ satisfies 
$X_{\mathrm{B}}[u_{1}](\tau) \subseteq \mathbb X_2$ for all $\tau \in [0,t]$. 
Now,~\eqref{eq::odeT}-\eqref{eq::initT} imply that $x \in \mathbb{X}[u_1]$ and, 
consequently, $x(t) \in X[u_1](t)$. Thus, we have established~\eqref{eq::inclusion1}.

For the second part of the proof, we need to show that
the inclusion 
\begin{equation}
\label{eq::inclusion2}
X[u_1](t) \subseteq X_{\mathrm{F}}[u_1](t) \; ,
\end{equation}
holds. Let $x \in \mathbb{X}[u_1]$ and $u_2$ satisfy
\begin{align}
\dot x(\tau) &= f(x(\tau),u_1(\tau),u_{2}(\tau)) \\
x(\tau) &\in \mathbb X_2 \\
x(0) &= 0 \; .
\end{align}
It is clear that $x$ satisfies the constraints in~\eqref{eq::breachset} on $[t,T]$, 
$x(t) \in X_{\mathrm{B}}[u_{1}](t)$. But then, $x$ also satisfies the 
constraints in~\eqref{eq::freachset} on $[0,t]$, which implies~\eqref{eq::inclusion2}.
Finally,~\eqref{eq::inclusion1} and~\eqref{eq::inclusion2} yield the statement of
the theorem.
\hfill\hfill$\diamond$

\subsection{Construction of convex enclosures using generalized differential 
inequalities}
\label{sec::gdi}

This section is concerned with the construction of enclosures for the
reachability tube $X_{\rm F}[u_1]$. 

\begin{defn}
Let $Z:\mathbb{R}\to\mathcal{P}(\mathbb{R}^{n})$ be a set-valued function. 
A set-valued function $Y:\mathbb{R}\to\mathbb{K}^{n}_{\rm C}$ is called an 
enclosure of $Z$ on $[0,T]$ if $Y(t)\supseteq Z(t)$ for all $t\in[0,T]$.
\end{defn}

In the following, we use the shorthand
\begin{equation}
\label{eq:gammaDef}
\Gamma(\nu_1,c,Y,Z) = \left\{ f(\xi,\nu_1,\nu_2) \ \middle| \ 
\begin{alignedat}{1}
c^\intercal\xi &= V[Y](c) \\
\xi &\in Y \cap \mathrm{int}(Z) \\
\nu_2 &\in \mathbb{U}_2
\end{alignedat}
 \right\}\;,
\end{equation}
which is defined for all $Z \in \mathbb{K}^{n_{x}}_{\rm C}$, 
$c \in \mathbb R^{n_x}$, and $\nu_1 \in \mathbb R^{n_u}$.

The next theorem provides a basis for the construction of convex 
enclosures of $X_{\rm F}[u_1]$. It exploits the reach-set splitting
structure of the backward-forward propagation scheme that has been introduced in
Section~\ref{sect::bfRsS}.

\begin{thm}
\label{thm::outerreach}
Let Assumptions~\ref{ass::rhs1} and~\ref{ass::admissible} be satisfied and let 
the Lebesgue integrable function $u_1:\mathbb{R}\to\mathbb{U}_2$ be given. Let 
$\mathbb X_2,Y_{\rm B}, Y_{\rm F}: [0,T] \to \mathbb K_{\rm C}^{n_x}$ be compact
set-valued functions such that $V[Y_{\rm B}(\cdot)](c)$ and 
$V[Y_{\rm F}(\cdot)](c)$ are, for all $c \in \mathbb{R}^{n_x}$,
Lipschitz continuous on $[0,T)$. If
$$Y_{\rm B}(t) \cap \mathrm{int}(\mathbb{X}_2(t)) \neq \varnothing \quad \text{and} \quad 
Y_{\rm F}(t) \cap \mathrm{int}(Y_{\rm B}(t)) = \varnothing$$
for all $t \in [0,T]$ and if the inequalities
\begin{alignat*}{2}
&\dot V[Y_{\rm B}(t)](c) &&\leq -V[-\Gamma(u_1(t), c,Y_{\rm B}(t), \mathbb X_2(t))](c)  \\
&\dot V[Y_{\rm F}(t)](c) &&\geq V[\Gamma(u_1(t), c,Y_{\rm F}(t), Y_{\rm B}(t) )](c) \\
&V[Y_{\rm F}(0)](c) &&\geq c^\intercal x_{0} \\
&V[Y_{\rm B}(T)](c) &&\geq V[\mathbb X_2(T)](c)
\end{alignat*}
hold for all $t \in [0,T)$ and all $c\in\mathbb{R}^{n_x}$, then
the set-valued function $Y_{{\rm F} \cap {\rm B}}:\mathbb{R}\to\mathbb{K}^{n_x}_{\rm C}$
given by 
\begin{equation*}
\forall t\in[0,T], \quad Y_{{\rm F} \cap {\rm B}}(t) = Y_{\rm F}(t) \cap Y_{\rm F}(t)
\end{equation*}
is an enclosure of $X[u_1]$ on $[0,T]$.
\end{thm}

\begin{figure*}[htbp!]
        \centering

\begin{overpic}[width=0.9\textwidth]{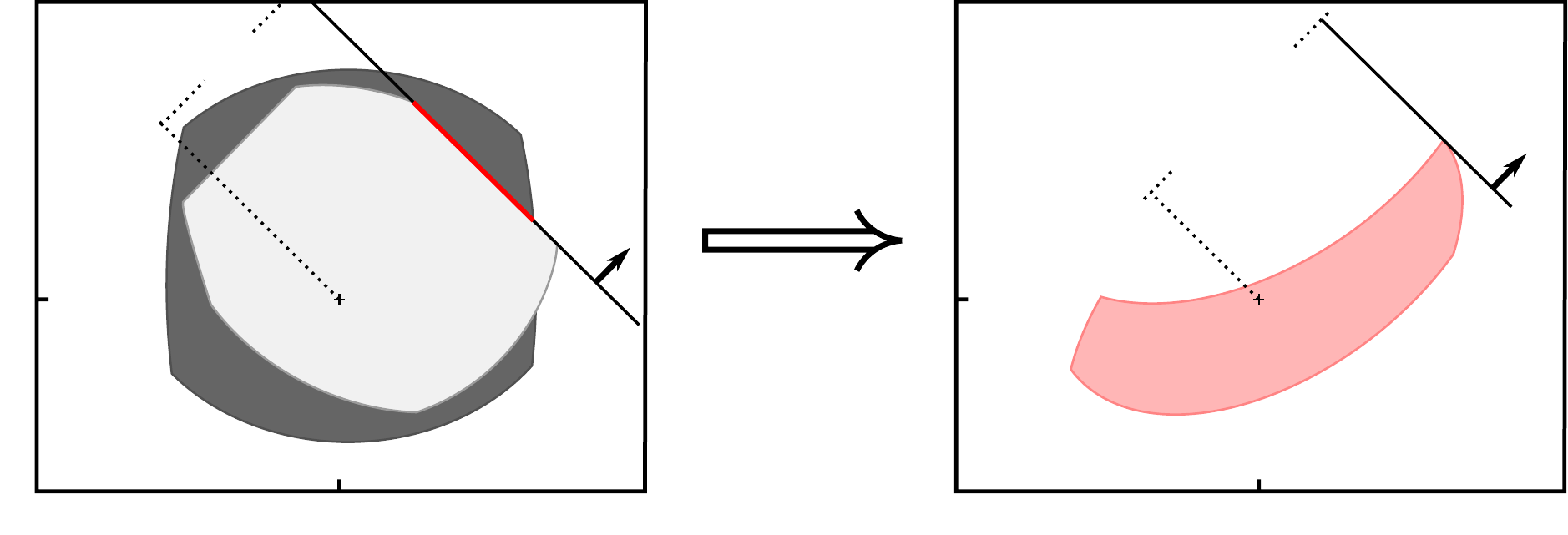}
\put(0.5,14.5){\footnotesize $0$}
\put(21.2,0.5){\footnotesize $0$}
\put(80,0.5){\footnotesize $0$}
\put(58.9,14.5){\footnotesize $0$}

\put(3.2,31.7){\footnotesize $x_2$}
\put(38.7,4){\footnotesize $x_1$}

\put(62,31.7){\footnotesize $f_2$}
\put(97.5,4){\footnotesize $f_1$}

\put(45,22){\footnotesize $f(\cdot,u_1(t),\mathbb{U}_2)$}

\put(8,30){\footnotesize $V[Y_{\rm F}(t)](c)$}
\put(65,26){\footnotesize $V[\Gamma(u_1(t),c,Y_{\rm F}(t),Y_{\rm B}(t))](c)$}

\end{overpic} 

\caption{\label{fig::Ysketch}
A sketch of the conditions in Theorem~\ref{thm::outerreach}.
Left: The sets $Y_{\rm F}(t)$ and $Y_{\rm B}(t)$ are shown in light and dark 
gray. The set $\{\xi \in Y_{\rm F}(t) \ |\ c^\intercal \xi = V[Y_{\rm F}(t)](c),\,
\xi \in \operatorname{int}(Y_{\rm B}(t))\}$, is shown in red. Right: The set 
$\Gamma(u_1(t),c,Y_{\rm F}(t),Y_{\rm B}(t))$ is shown in light red. In both plots, 
the black arrow is the direction vector $c$. We also use the shorthand 
$f(\xi,u_{1}(t),\mathbb{U}_2) = \{ f(\xi,u_1(t),\nu_2) \ | \ \nu_2\in \mathbb{U}_2 \} $.}
\end{figure*}
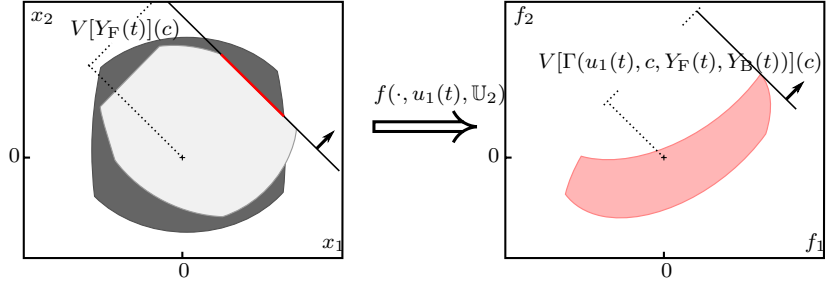

Figure~\ref{fig::Ysketch} visualizes the conditions of Theorem~\ref{thm::outerreach}. The left panel shows two dimensional enclosures $Y_{\rm F}(t)$ and $Y_{\rm B}(t)$ in light and dark gray, respectively. Notice that the conditions in Theorem~\ref{thm::outerreach} 
are verified pointwise-in-time at the boundary of the enclosures. The black arrow indicates a given direction $c$ and the red line corresponds to the set of points of the associated supporting facet,
\begin{equation*}
F[Y_{\rm F}(t)](c) = \left\{\xi  \ \middle|\ 
\begin{aligned}
c^\intercal \xi &= V[Y_{\rm F}(t)](c) \\
\xi &\in Y_{\rm F}(t)  
\end{aligned}\right\} \;,
\end{equation*}
that are also in $\operatorname{Y_{\rm B}(t)}$. For a given control $u_1$ and
starting from $F[Y_{\rm F}(t)](c) \cap \operatorname{int}(Y_{\rm B}(t))$, 
we have that any trajectory $x$ must satisfy 
\begin{equation*}
\begin{aligned}
\dot{x}(t) &\in \bigcup_{\xi\in F[Y_{\rm B}(t)](c) \cap \operatorname{int}(Y_{\rm B}(t))} \{ f\left(\xi,u_1(t),\nu_2\right) \ | \ \nu_2\in \mathbb{U}_2 \} \;.
\end{aligned}
\end{equation*}
At this point, it is easy to see that the right-hand side of the above inclusion 
is exactly $\Gamma(u_1(t),c,Y_{\rm F}(t),Y_{\rm B}(t))$---which is shown in the 
right panel, in light red. Then, taking the support function of
$\Gamma(u_1(t),c,Y_{\rm F}(t),Y_{\rm B}(t))$ as a bound on $\dot{V}[Y_{\rm F}(t)](c)$
we are bounding the dynamics of the points on 
$F[Y_{\rm F}(t)](c) \cap \operatorname{int}(Y_{\rm B}(t))$.

A proof of Theorem~\ref{thm::outerreach} can be found in Appendix~\ref{app::igdi}.

\subsection{Ellipsoidal-valued enclosures for reachability tubes}
\label{sec::ellgdi}
This section presents a practical construction of ellipsoidal enclosures
based on Theorem~\ref{thm::outerreach}. Our focus is on ellipsoidal
set parameterizations of the form
\begin{equation*}
Y_{\rm B}(t) = \mathcal{E}(q_{\rm B}(t),Q_{\rm B}(t)) \quad\text{and}\quad 
Y_{\rm F}(t) = \mathcal{E}(q_{\rm F}(t),Q_{\rm F}(t))\;.
\end{equation*}
In particular, our goal is to develop a computational method for
constructing the functions 
$q_{\rm B},q_{\rm F}:\mathbb{R}\to\mathbb{R}^{n_x}$ as well as
$Q_{\rm B},Q_{\rm F}:\mathbb{R}\to\mathbb{S}^{n_x}_{++}$
in such a way that $Y_{\rm B}$ and $Y_{\rm F}$ satisfy the conditions from 
Theorem~\ref{thm::outerreach}. 

\begin{assum}
\label{ass::closed}
The sets $\mathbb X_2(t)$ are bounded for all $t\in[0,T]$.
\end{assum}

Let $s:\mathbb{R}\to\mathbb{R}^{n_x}$, $S:\mathbb{R}\to\mathbb{S}^{n_x}_{++}$, and
$(v,V)\in\mathbb{R}^{n_u}\times \mathbb{S}^{n_u}_{++}$ be given such that
\begin{equation*}
\mathbb{U}_2 \subseteq \mathcal{E}\left( v, V \right) \quad \text{and} \quad 
\mathbb{X}_2(t) \subseteq \mathcal{E}\left( s(t), S(t) \right) \; .
\end{equation*}
The existence of $s,S,v$ and $V$ is guaranteed, if 
Assumptions~\ref{ass::admissible} and~\ref{ass::closed} hold.
Moreover, let 
\begin{equation*}
\Omega:\mathbb{R}^{n_x\times n_x} 
\times \mathbb{R}^{n_x\times n_u} \times \mathbb{R}^{n_x}\times\mathbb{R}^{n_u} 
\times\mathbb{R}^{n_u} \times \mathbb{S}^{n_x}_{+} \to \mathbb{S}^{n_x}_{+}
\end{equation*}
be a nonlinearity bound such that 
\begin{equation*}
f(x,u_1,u_2) - A (x-q) - B(u_2-v) \in \mathcal E(0,\Omega(A,B,q,u_1,v,Q))
\end{equation*}
is satisfied for all vectors $x \in \mathcal{E}(q,Q)$; all vectors $u_1$, $u_2$, 
$v$, and $q$; and all matrices $A$, $B$, and $Q$ of compatible size. If Assumption~\ref{ass::rhs1} is satisfied, such a function can always be 
constructed~\cite{Houska2012}.

\begin{rem}
\label{rem::nonlinearity}
The accuracy of the enclosures constructed in this section 
depends on the choice of $v$, $V$, $s$, $S$, and $\Omega$. 
A thorough analysis of the conservatism of the ellipsoidal bounds as a function 
of these parameters goes beyond the scope of this paper. Methods to construct 
such functions can be found in other works~\cite{Houska2012,Villanueva2017a}.
\end{rem}

In the following, we introduce the functions
\begin{equation}
\begin{alignedat}{1}
\Phi_{1}(Q,A) &= A Q + Q A^\intercal  \notag \\
\Phi_{2}(Q,W,\sigma) &= \sigma Q + \frac{1}{\sigma} W  \notag \\
\varphi_{3}(q_1,q_2,Q_1,Q_2,\kappa) &= \kappa Q_1Q_2^{-1}(q_2-q_1) \notag \\
\Phi_{3}(q_1,q_2,Q_1,Q_2,\kappa) &= \kappa \left( 
I - \left\Vert q_1-q_2 \right\Vert_{Q_2^{-1}}^{2} I - Q_1 Q_2^{-1} \right) Q_1  \notag
\end{alignedat}
\end{equation}
for all scalars $\kappa,\sigma$ as well as all vectors $q_1,q_2$ and matrices 
$A,B,W,Q,Q_1,Q_2$ with compatible dimensions. Similarly, we introduce the 
variables 
\begin{equation*}
y = ( q_{\rm B}, q_{\rm F}, Q_{\rm B}, Q_{\rm F} ) \in \mathbb Y
\end{equation*}
and
\begin{equation*}
\lambda = \left( A_{\rm F}, A_{\rm B}, B_{\rm F}, B_{\rm B}, \sigma_{\rm B}, 
\sigma_{\rm F}, \mu_{\rm B}, \mu_{\rm F}, \kappa_{\rm B}, \kappa_{\rm F} \right)
\in \mathbb L 
\end{equation*}
together with the domains $\mathbb{Y}= \mathbb{R}^{n_x} \times \mathbb{R}^{n_x} 
\times \mathbb{S}^{n_x}_{++} \times \mathbb{S}^{n_x}_{++}$ and 
$\mathbb{L}=\mathbb{R}^{n_x\times n_x} \times \mathbb{R}^{n_x\times n_x} \times 
\mathbb{R}^{n_x \times n_u} \times \mathbb{R}^{n_x \times n_u} \times \mathbb{R}^{6}_{+}$.
The functions $F:\mathbb{Y} \times \mathbb{R}^{n_u} \times \mathbb{R}^{n_x} \times \mathbb{S}^{n_x}_{++} \times \mathbb{L} \to \mathbb{Y}$ with
$F = (F_1,F_2,F_3,F_4)$ and 
\begin{alignat*}{2}
&F_1( t,y,u_1,\lambda ) &&= \hphantom{{}+{}} f( q_{\rm B},u_1,v ) 
- \varphi_{3}(q_{\rm B},s(t),Q_{\rm B},S(t),\kappa_{\rm B}) \\
&F_2( t,y,u_1,\lambda ) &&= \hphantom{{}+{}} f( q_{\rm F},u_1,v ) 
+ \varphi_{3}(q_{\rm F},q_{\rm B},Q_{\rm F},Q_{\rm B},\kappa_{\rm F}) \\
&F_3(t,y,u_1,\lambda) &&= \hphantom{{}+{}} \Phi_{1}(Q_{\rm B},A_{\rm B}) 
- \Phi_{2}\left(Q_{\rm B},BVB^\intercal,\sigma_{\rm B},\right) \\
& &&\hphantom{{}={}} - \Phi_{2}\left(Q_{\rm B},\Omega(A_{\rm B},B_{\rm B},q_{\rm B},u_1,v,Q_{\rm B}),
\mu_{B}\right) \\
& &&\hphantom{{}={}} - \Phi_{3}\left(q_{\rm B},s(t),Q_{\rm B},S(t),\kappa_{\rm B}\right) \\
&F_4(t,y,u_1,\lambda) &&= \hphantom{{}+{}} \Phi_{1}(Q_{\rm F},A_{\rm F}) 
+ \Phi_{2}\left(Q_{\rm F},BVB^\intercal,\sigma_{\rm F}\right) \\
& &&\hphantom{{}={}} + \Phi_{2}\left(Q_{\rm F},\Omega(A_{\rm F},B_{\rm F},q_{\rm F},u_1,v,Q_{\rm F}),
\mu_{\rm F}\right) \\
& &&\hphantom{{}={}} + \Phi_{3}\left(q_{\rm F},q_{\rm B},Q_{\rm F},Q_{\rm B},
\kappa_{\rm F}\right)
\end{alignat*}
and
$G:\mathbb{Y} \times \mathbb{Y} \to\mathbb{Y}$ given by
\begin{equation*}
\begin{aligned}
&G( y(0), y(T) ) \\
&\quad = 
\left( q_{\rm F}(0) - x_0 , \, q_{\rm B}(T) - s(T)  , \, 
Q_{\rm F}(0), \, Q_{\rm B}(T) - S(T) \right) \;.
\end{aligned}
\end{equation*}
determine, respectively, the right-hand side and constraint function
of the boundary-value problem that is needed in the following theorem.

\begin{thm}
\label{thm::ellreach}
Let Assumptions~\ref{ass::rhs1},~\ref{ass::admissible}, and~\ref{ass::closed} be
satisfied and let the Lebesgue integrable function $u_1: [0,T] \to \mathbb U_1$ be given. 
Let $y:\mathbb{R}\to\mathbb{Y}$ and 
$\lambda: \mathbb{R} \to \mathbb L$ be any functions satisfying the boundary
value problem
\begin{gather*}
\begin{alignedat}{3}
&\forall t\in[0,T], \quad 
&&\dot{y}(t) &&= F\left( t, y(t), u_1(t), \lambda(t) \right) \\
& &&0 &&= G( y(0), y(T) ) \; .
\end{alignedat}
\end{gather*}
Then, the set-valued functions
$$Y_{\rm B}(t) = \mathcal{E}(q_{\rm B}(t),Q_{\rm B}(t)) \quad \text{and} \quad Y_{\rm F}(t) = \mathcal{E}(q_{\rm F}(t),Q_{\rm F}(t))$$
satisfy the conditions of Theorem~\ref{thm::outerreach} on $[0,T]$. That is, the 
set-valued function $Y_{{\rm F}\cap{\rm B}}:\mathbb{R}\to\mathbb{K}^{n_x}_{\rm C}$ with 
\begin{equation*}
\forall t\in[0,T], \quad Y_{{\rm F}\cap{\rm B}}(t) =  
\mathcal{E}(q_{\rm F}(t),Q_{\rm F}(t))\cap\mathcal{E}(q_{\rm B}(t),Q_{\rm B}(t))
\end{equation*}
 is an enclosure of $X[u_1]$ on $[0,T]$.
\end{thm}

A proof of Theorem~\ref{thm::ellreach} can be found in Appendix~\ref{app::ellreach}.

\section{Tractable approximation of differential games}
\label{sec::case}

We now present a conservative and tractable approximation of~\eqref{eq::exactgame}.
This approximation is constructed by leveraging on the properties of the ellipsoidal
approximation presented in the previous section. We assume that two bounding functions
$H: \mathbb R \times \mathbb{R}^{n_y} \to \mathbb{R}^{n_h}$ and $M: \mathbb{R}^{n_y}
\to \mathbb{R}^{n_h}$ satisfying 
\begin{align}
H(t,y) \leq 0 \; \; &\Longrightarrow \; \; \mathcal E( q_{\rm F}, Q_{\rm F} ) \cap \mathcal E( q_{\rm B}, Q_{\rm B} ) \subseteq \mathbb{X}_{1}(t) \\
M(y) \leq 0 \; \; &\Longrightarrow \; \; \mathcal M( \mathcal E( q_{\rm F}, Q_{\rm F} ) \cap \mathcal E( q_{\rm B}, Q_{\rm B} ) ) \leq 0 \; ,
\end{align}
are available. Notice that the construction of these functions is akin 
to the construction of the nonlinearity bounder $\Omega$~\cite{Villanueva2018}; see also Remark~\ref{rem::nonlinearity}.

Let us consider the optimal control problem
\begin{equation}
\label{eq::ellgame}
\begin{alignedat}{2}
&\inf_{ \substack{ x, y, \\ u_1, u_2, \lambda} } \ M( y(T) ) \quad  \operatorname{s.t.} \,
\left\lbrace
\begin{array}{l}
\forall t \in [0,T], \\
\dot{x}(t)= f( x(t),u_1(t),u_2(t) ) \\
\dot{y}(t) = F( y(t),u_1(t),\lambda(t) ) \\ 
u_{1}(t) \in \mathbb{U}_1 \, ,\ u_{2}(t) \in \mathbb{U}_2 \\
y(t) \in \mathbb{Y} \, ,\ \lambda(t) \in \mathbb{L} \\
x(t)\in\mathbb{X}_2(t) \\
0 \geq H( t, y(t) ) \\
0 = G(y(0),y(T))  \;.
\end{array}
\right . 
\end{alignedat}
\end{equation}
Theorem~\ref{thm::ellreach} implies that any feasible point of
is a feasible point of~\eqref{eq::exactgame}. The auxiliary state $x$ is used 
to enforce lower-level feasibility of the control $u_1$.

\subsection{Numerical illustration}
\label{sec::numerics}

We consider a differential game with dynamics given by 
\begin{equation*}
\begin{alignedat}{2}
&\forall t\in[0,T],\quad
&&\dot{x}(t) = 
\begin{pmatrix}
&x_1(t) + \frac{1}{2} x_2(t) + u_1(t),\\
&\frac{3}{2} x_1(t) + x_2(t) + u_2(t)
\end{pmatrix}\\[0.1cm]
&\text{and} &&x(0) = (0,0)^\intercal
\end{alignedat}
\end{equation*}
with $T = \frac{3}{2}$.
The path constraints are given by
\begin{equation*}
\mathbb X_1(t) = [-6, 6]^{2} \quad \text{and} \quad \mathbb X_2(t) = R(t)\mathcal E(s, S)
\end{equation*}
for all $t \in [0,T]$ with
\begin{equation*}
\begin{gathered}
s = 
\begin{pmatrix}
\frac{5}{50}\\[0.2cm]
-\frac{71}{25}
\end{pmatrix}\;, \quad S = 
\begin{pmatrix}
\frac{802}{25} & \frac{16}{5}\\[0.1cm]
\frac{16}{5} & \frac{802}{25}
\end{pmatrix}\;, \\[0.1 cm]
\text{and}\quad R(t) =
\begin{pmatrix}
\frac{21}{20}\cos(\frac{\pi t}{2}) & -\frac{19}{20}\sin(\frac{\pi t}{2})\\[0.2cm]
\frac{21}{20}\sin(\frac{\pi t}{2}) & \frac{19}{20}\cos(\frac{\pi t}{2})
\end{pmatrix} \;.
\end{gathered}
\end{equation*}
The control sets are $\mathbb U_1 = [-1.4, 0]$ and $\mathbb U_2 = [0, 2]$. The 
Mayer objective function is given by $m(x) = x_{1} -x_{2} $\;.

The bounding functions $M$ and $H$ can, in this example, be constructed without introducing further conservatism using the result from Proposition~\ref{prop::intersection} in Appendix~\ref{app::ellreach}.
Problem~\eqref{eq::ellgame} was formulated and solved numerically with 
\texttt{ACADO Toolkit}~\cite{Houska2011} using a 
multiple shooting discretization with $10$ equidistant intervals, 
and a Runge-Kutta integrator of order 4/5.

\begin{figure}[htbp!]
        \centering

\begin{overpic}[width=0.4\columnwidth]{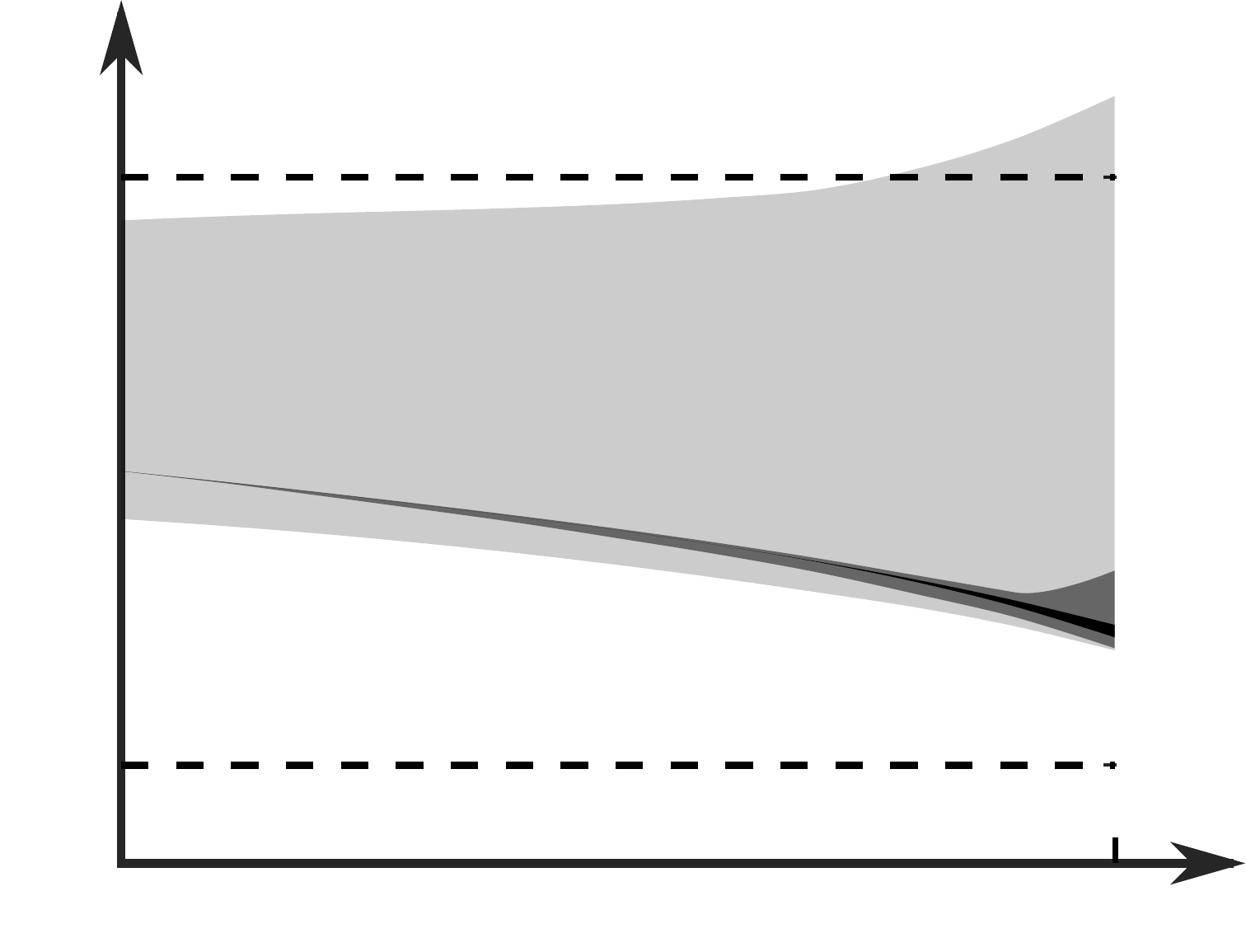}
\put(8,2){\footnotesize $0$}
\put(88,2){\footnotesize $T$}
\put(2,13.5){\footnotesize $-6$}
\put(5,61){\footnotesize $6$}
\put(50,-3){\small $t$}
\put(-3,40){\small $x_1$}
\end{overpic} 

\vspace{0.5cm}
 
\begin{overpic}[width=0.4\columnwidth]{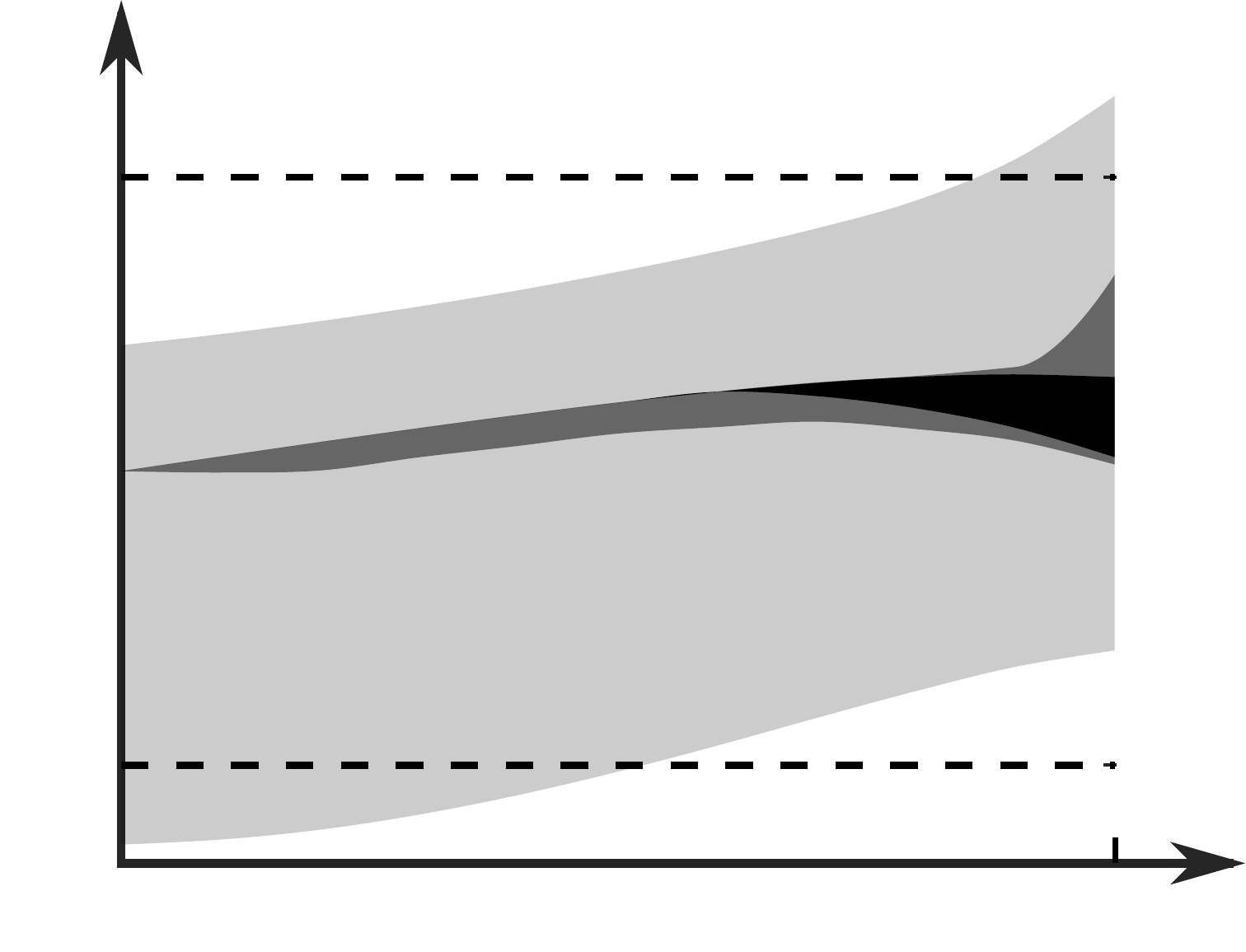}
\put(8,2){\footnotesize $0$}
\put(88,2){\footnotesize $T$}
\put(2,13.5){\footnotesize $-6$}
\put(5,61){\footnotesize $6$}
\put(50,-3){\small $t$}
\put(-3,40){\small $x_2$}
\end{overpic} 

\vspace{0.5cm}

\caption{\label{fig::tube} Projections of $Y_{\rm B}$ (light gray) and 
the pointwise-in-time intersection of $Y_{\rm B}$ and $Y_{\rm F}$ (dark gray)
onto the $x_1$- and $x_2$-axis.}
\end{figure}
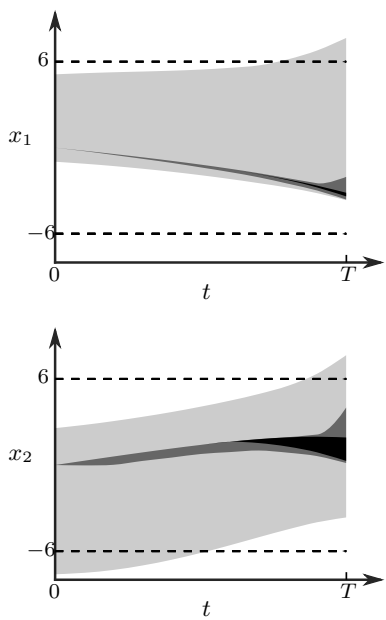

Figure~\ref{fig::tube} shows projections onto the $x_1$- (top) and $x_2$-axis (bottom)
planes of the ellipsoidal enclosures $Y_{\rm B}$ (light gray), the pointwise-in-time
intersection of $Y_{\rm B}$ and $Y_{\rm F}$ (dark gray), and an inner
approximation of $X[u_1^{\star}]$ (black)---computed by Monte Carlo 
simulation with $10^4$ trajectories.

The optimal value of $(26)$ is $-1.59$.
Recall that~\eqref{eq::ellgame} is only a conservative approximation 
of~\eqref{eq::exactgame}, since Theorem~\ref{thm::ellreach} only provides a means
to construct outer approximations of the exact forward reachable set of the game. 
However, the conservatism of the solution can be evaluated a posteriori by solving the
optimal control problem of Player~2. Here, we find that Player~2
can at most achieve an optimal value of $-3.67$. Thus, the optimal value of the exact differential game is overestimated by approximately~$2.08$. If one wishes to further reduce this overestimation, one would have to abandon the idea to work with ellipsoids and use more accurate set parameterizations. An in depth analysis of such general set parameterizations is, however, beyond the scope of this paper.

Last but not least, we also solve the above game for the case that Player~2 has no state constraints. In this case, the differential game reduces to a standard robust optimal control problem and the optimal value of Player~2 is $3.85$.
This illustrates the importance of taking state constraints in dynamic games into account and highlights the differences between robust optimal control and more general
zero-sum differential games with state constraints.

\section{Conclusion}
\label{sec::conclusion}

This paper has presented a set-theoretic framework for the numerical analysis of
zero-sum differential games with state constraints. In particular, it introduced 
a novel backward-forward reachable set splitting scheme, which can be used by 
the first player to compute the reachable set of states of the second player. 
This splitting scheme was then used to derive convex outer approximations for the
reachable set of the game using boundary value constrained differential inequalities.
A particular emphasis was placed on ellipsoidal outer approximations, 
which lead to the conservative but tractable approximation~\eqref{eq::ellgame}
of the solution of the original Stackelberg differential game~\eqref{eq::exactgame}.
The advantage of~\eqref{eq::ellgame}
is that this is a standard optimal control problem, which can be solved using
state-of-the art optimal control solvers.
The effectiveness of the approach was illustrated by means of a numerical
example.

\section*{Acknowledgements}                               
\noindent
This research was supported by the National Natural Science Foundation China, Grant-No.~61473185; as well as ShanghaiTech University,
Grant-No.~F-0203-14-012.

\bibliographystyle{plain}
\bibliography{games}

\appendix
\section{\label{app::igdi}Proof of Theorem~\ref{thm::outerreach}}

The proof of Theorem~\ref{thm::outerreach} is non-trivial and, therefore, it has been divided into different sub-sections, which build upon each other.

\subsection{Technical preliminaries}
\label{subsec::prelim}
Let us consider a general differential equation of the form
\[
\dot x(t) = g(x(t), w(t)) \quad \text{with} \quad x(0) = x_0
\]
with $g$ being Lipschitz continuous in $x$ and continuous in the external input 
$w: [0,T] \to \mathbb W \subseteq \mathbb{K}^{n_w}$. Let
\begin{equation*}
X_{\mathbb T}(t) = \left\{
\xi \in \mathbb R^{n_x}
\middle |
\begin{aligned}
&\exists x \in \mathbb W^{n_x}_{1,1},\ \exists w \in \mathbb T^{n_w}: \\
&\aew \tau \in [0,t], \\
&\dot x(\tau) = g(x(\tau),w(\tau))\\
& x(0) = x_0 \\
&w(\tau) \in \mathbb W, \ x(t)  = \xi
\end{aligned}
\right\}
\end{equation*}
denote the reachable sets of this differential equation for a given subset 
$\mathbb T \subseteq \mathbb L_{1}$. Let
$\mathbb C^{n_w} \subseteq \mathbb L_1^{n_w}$  denote the set of bounded continuous
functions. Now, a direct consequence of Lusin's theorem~\cite{Feldman1981} and 
Gronwall's lemma is that
\begin{equation*}
\forall t \in \mathbb R, \qquad \text{cl}(X_{\mathbb C}(t)) = \text{cl}(X_{\mathbb L_1}(t)) \; .
\end{equation*}
In other words, if we are only interested in the closure of a reachable set, we may
simply replace Lebesgue integrable functions by continuous functions in a statement
without altering its conclusion. In the following, we will use such replacements
without saying this explicitly at all places. In particular, we assume, without 
loss of generality, that $u_1$ is any given continuous function.

\subsection{Constrained set propagation operators}
Let
\begin{equation}
\label{eq::PI}
\Pi( t_2,t_1, X_1, Z ) = \left\{
\xi \in \mathbb R^{n_x}
\middle |
\begin{aligned}
&\exists x \in \mathbb W^{n_x}_{1,1},\ \exists u_2 \in \mathbb L^{n_u}_{1}: \\
&\aew t \in [t_1,t_2], \\
&\dot x(t) = f(x(t),u_1(t),u_2(t))\\
&x(t_1) \in X_1, \ x(t) \in Z(t), \\
&u_2(t) \in \mathbb U_2, \ x(t_2)  = \xi
\end{aligned}
\right\} \; ,
\end{equation}
denote the constrained set propagation operator of~\eqref{eq::odes}, which is
defined for all $X_{1}\in\mathbb{K}_{\rm C}^{n_x}$, all 
$Z: [t_1,t_2] \to \mathbb{K}_{\rm C}^{n_x}$,
and all $t_1,t_2\in\mathbb{R}$ with $t_1\leq t_2$.
Moreover, let
\begin{equation*}
\pi\left( \xi, X \right) = \operatorname*{argmin}_{\xi^{'} \in X} 
\left\Vert \xi - \xi^{'} \right\Vert_2^2 \; .
\end{equation*}
denote the Euclidean projection of a point $\xi \in \mathbb R^{n_x}$ onto a 
compact convex set $X \in \mathbb{K}_{\rm C}^{n_x}$. Now, the key idea is to 
introduce the auxiliary differential equation
\begin{align}
\label{eq::ODE2}
\forall t \in [0,T], \quad \dot z(t) = f_{K}( t, z(t), u_2(t), Z(t) )
\end{align}
with
\begin{equation*}
f_{K}( t, \xi, \nu_2, X ) = - K ( \xi - \pi( \xi, X ) ) + f( \xi, u_1(t), \nu_2 )
\end{equation*}
for all $\xi\in\mathbb{R}^{n_x}$, $\nu_2\in\mathbb{R}^{n_u}$, and all 
$X\in\mathbb{K}_{\rm C}^{n_x}$.
Here, $K>0$ is a tuning parameter that can be interpreted as a proportional control
gain of an additional control term, which can be used to steer~\eqref{eq::ODE2} towards 
$X$ whenever $z(t)$ is outside of $X$. Let
\begin{equation*}
\Pi_{K}(t_2,t_1,X_1,Z) = \left\{
\xi \in \mathbb R^{n_x}
\middle |
\begin{aligned}
&\exists z \in \mathbb W^{n_x}_{1,1},\ \exists u_2 \in \mathbb L^{n_u}_{1}: \\
&\aew t \in [t_1,t_2], \\
&\dot z(t) = f_{K}(t,z(t),u_2(t),Z(t))\\
&x(t_1) \in X_1, \ u_2(t) \in \mathbb U_2\\ 
&x(t_2)  = \xi
\end{aligned}
\right\}\;,
\end{equation*}
denote the set-propagation operator of~\eqref{eq::ODE2}. In analogy to the 
propagation operator $\Pi$, $\Pi_K$ is defined for all 
$X_{1}\in\mathbb{K}^{n_x}_{\rm C}$, all $Z: [t_1,t_2] \to \mathbb{K}_{\rm C}^{n_x}$,
and all $t_1,t_2\in\mathbb{R}$ with $t_1\leq t_2$.

\begin{lem}
\label{lem::GGDI}
Let Assumptions~\ref{ass::rhs1} and~\ref{ass::admissible} be satisfied.
Let $Y,Z:~[0,T] \to \mathbb{K}_{\rm C}^{n_x}$ be any given set valued function
such that the intersection $Y(t)\cap Z(t)$ is, for all $t\in[0,T]$, nonempty
and such that the functions $V[Y(\cdot)\cap Z(\cdot)](c)$ and $V[Y(\cdot)](c)$ 
are, for all $c\in\mathbb{R}^{n_x}$, Lipschitz continuous on $[0,T]$. If the 
differential inequality
\begin{align}
\label{eq::AUXgdi}
&\begin{alignedat}{2}
\dot V[Y(t)](c) \geq & \max_{\xi,\nu_2} \ 
&& c^\intercal f_{K}(t,\xi,\nu_2, Y(t) \cap Z(t) ) \\
&  \operatorname{s.t.} && \left\{
\begin{aligned}
c^\intercal \xi &= V[Y(t)](c) \\
\xi &\in Y(t) \\
\nu_2 &\in \mathbb U_2 \\
\end{aligned} \right.
\end{alignedat}  \\[0.1cm]
&V[Y(0)](c) \geq V[X_0](c)
\end{align}
is satisfied for all $c \in \mathbb R^{n_x}$ and all $t \in [0,T]$
for a given initial set $X_{0}\in\mathbb{K}^{n_x}_{\rm C}$, then
\begin{equation*}
\forall t\in[0,T],\quad \Pi(0, t, X_0, Z ) \subseteq  \Pi_{K}(0,t,X_0,Z)
\subseteq Y(t) \;.
\end{equation*}
\end{lem}

\textbf{Proof.}
First, notice that $\pi$ satisfies $\pi(\xi,Y(t)\cap Z(t)) = \xi$ whenever 
$\xi\in Y(t)\cap Z(t)$. Thus, we have
\begin{equation*}
f_{K}(t,\xi,u_2(t), Y(t) \cap Z(t)) = f(\xi,u_1(t), u_2(t))
\end{equation*}
for all $\xi\in Y(t)\cap Z(t)$. This implies the first inclusion,
\begin{equation*}
\forall t\in[0,T], \qquad \Pi(0, t, X_0, Z ) \subseteq \Pi_{K}(0,t,X_0,Z) \; .
\end{equation*}
In order to establish the remaining inclusion, we assume for a moment that
$u_2$ is constant.
Since the set $Y(t)\cap Z(t)$ is nonempty and convex, $\pi(\cdot,Y(t)\cap Z(t))$ 
is non-expansive. Thus, this function is uniformly Lipschitz continuous on $[0,T] \times \mathbb{U}_2$ with Lipschitz constant~$1$. 
Furthermore, by the Lipschitz continuity of the function $V[Y(\cdot)\cap Z(\cdot)](c)$ the function $\pi(\xi,Y(\cdot)\cap Z(\cdot))$ is continuous. This, together with
Assumption~\ref{ass::rhs1}, implies that $f_{K}$ is jointly continuous in $(t,\xi,\nu_2)$ as well as Lipschitz continuous in $\xi$, uniformly on $[0,T]\times\mathbb{U}_2$. 

At this point, it is important to notice that $Y$ and $Z$ are arbitrary but given, thus 
the right-hand side function $f_{K}$ is only a function of time, the state and the parameter $u_2$. We can now replicate the arguments of all the steps in the proof of
Theorem~3 in~\cite{Villanueva2015} to~\eqref{eq::ODE2}---observing that we have 
established the required properties of $f_K$---to obtain
the second inclusion 
\begin{equation*}
\forall t\in[0,T],\quad \Pi_{K}(0,t,X_{0},Z) \subseteq Y(t).
\end{equation*}
The assumption that $u_2$ is constant can be removed using the same argument
as in~\cite[Remark~2]{Villanueva2015}.\hfill\hfill$\diamond$

\begin{cor}
\label{cor::GGDI1}
Let the conditions of Lemma~\ref{lem::GGDI} hold. Assume, in addition, that 
$Y(t)\cap\operatorname{int}(Z(t))\neq \varnothing$, for all $t\in[0,T]$. If $Y$ 
is such that the differential inequality
\begin{align*}
&\begin{alignedat}{2}
\dot V[Y(t)](c) \geq & \max_{\xi,\nu_2} \ && c^\intercal f(\xi,u_1(t),\nu_2) \\
&  \operatorname{s.t.} && \left\{
\begin{aligned}
c^\intercal \xi &= V[Y(t)](c) \\
\xi &\in Y(t) \cap \operatorname{int}(Z(t)) \\
\nu_2 &\in \mathbb U_2 \\
\end{aligned} \right.
\end{alignedat}  \\[0.1cm]
&V[Y(0)](c) \geq V[X_0](c)
\end{align*}
is satisfied for all $c \in \mathbb R^{n_x}$ and all $t \in [0,T]$,
then
$$\forall t \in [0,T], \qquad \Pi(0,t,X_0,Z) \subseteq Y(t) \; .$$
\end{cor}

\textbf{Proof.}
The proof proceeds in two steps. First, we show that the result
holds under the stronger assumption that $Y$ is such that $Y(t)$ is strictly 
convex for all $t\in[0,T]$ and that it satisfies the differential inequality
\begin{equation}
\label{eq::GDIaux}
\begin{alignedat}{2}
\dot V[Y(t)](c) \geq & \max_{\xi,\nu_2} \ && c^\intercal f(\xi,u_1(t),\nu_2) \\
&  \operatorname{s.t.} && \left\{
\begin{aligned}
c^\intercal \xi &= V[Y(t)](c) \\
\xi &\in Y(t) \cap Z(t) \\
\nu_2 &\in \mathbb U_2 \\
\end{aligned} \right.
\end{alignedat}  
\end{equation}
for all $t\in[0,T]$ and all $c\in\mathbb{R}^{n_x}$. 

Let the set 
\begin{equation*}
F[Y(t)](c) = \left\{ \xi\in\mathbb{R}^{n_x} \middle | 
\begin{aligned}
c^\intercal \xi & = V[Y(t)](c) \\
\xi &\in Y(t)
\end{aligned}
 \right\} \;.
\end{equation*}
be the supporting facet of $Y(t)$ in the direction 
$c\in\mathbb{R}^{n_x}$. Notice that for any given pair 
$(t,c)\in[0,T]\times\mathbb{R}^{n_x}$, the set $F[Y(t)](c)$ is a singleton. Thus, for any given pair $(t,c) \in [0,T ] \times \mathbb{R}^{n_x}$, there are only two possible cases:

\textit{Case 1:} The set $F[Y(t)](c) \cap Z(t)$ is nonempty.
In this case, we have 
\begin{equation*}
F[Y(t)](c) \cap Z(t) \subseteq Y(t) \cap Z(t)\;.
\end{equation*}
Since $f$ and $f_{K}$ coincide on $Y(t) \cap Z(t)$ they also coincide 
on $F[Y(t)](c) \cap Z(t)$. Thus, $Y$ and $Z$ satisfy the differential 
inequality~\eqref{eq::AUXgdi} from Lemma~\ref{lem::GGDI} at $(t,c)$ for any
$K>0$.

\textit{Case 2:} The set $F[Y(t)](c) \cap Z(t)$ is empty.
This is only possible, if
\begin{equation*}
\max_{\xi \in F[Y(t)](c)} \ c^\intercal \xi > 
\max_{ \xi^{'}\in Y(t) \cap Z(t) } c^\intercal \xi^{'}\;. 
\end{equation*}
Thus, it follows that
\begin{equation*}
\forall \xi \in F[Y(t)](c),\quad \ c^\intercal ( \xi - \pi(\xi,Y(t) \cap Z(t)) > 0 \;,  
\end{equation*} 
since $\pi(\xi,Y(t) \cap Z(t)) \in Y(t) \cap Z(t)$\;. Therefore, the term
\begin{align*}
& c^\intercal f_{K}(t,\xi,u_2(t),Y(t) \cap Z(t)) \\
&= -K \underbrace{c^\intercal (\xi - \pi(\xi,Y(t) \cap Z(t)))}_{>0} + c^\intercal f(\xi,u_1(t),u_2(t))
\end{align*}
can be made arbitrarily small by choosing a sufficiently large $K$. 

Thus, we have shown that there exists for every pair $(t,c)\in[0,T]\times\mathbb{R}^{n_x}$ a sufficiently large $K$, such that $Y$ and $Z$ satisfy
the differential inequality~\eqref{eq::AUXgdi} from Lemma~\ref{lem::GGDI}. In particular
$Y$ and $Z$ satisfy the (strengthened) conditions from Lemma~\ref{lem::GGDI} in the limit 
as $K\to \infty$. We must mention that one should be careful when taking this limit, as the Lipschitz constant of $f_{K}$ diverges for $K\to\infty$. Fortunately, one 
can apply the following topological argument: if $Y$ and $Z$ satisfy the above hypothesis, 
one can always construct enclosures $Y_{\epsilon}$ and $Z_{\epsilon}$ satisfying the 
conditions of Lemma~\ref{lem::GGDI} for a sufficiently large $K$ and such that the 
Hausdorff distance between $Y_{\epsilon}(t)$ and $Y(t)$ (as well as $Z_{\epsilon}(t)$ and
$Z(t)$) converges to zero as $\epsilon\to 0$---uniformly on $[0,T]$. This claim follows 
readily from our continuity assumptions. As the images of these functions are compact, 
one can pass to the topological closure to show that Lemma~\ref{lem::GGDI} implies that
\begin{equation*}
\forall t \in [0,T], \qquad \Pi(0,t,X_0,Z) \subseteq Y(t) \; ,
\end{equation*}
if~\eqref{eq::GDIaux} holds.

At this point, we construct an enclosure 
$Y_{\epsilon}$ of the operator $\Pi$, such that 
$Y_{\epsilon}(t)$ is, for all $t\in[0,T]$, compact, and strictly convex. Then, 
we apply the procedure above and a continuity argument to show 
that~\eqref{eq::GDIaux} also holds in the limit as $\epsilon \to 0$, for 
set-valued functions with convex and compact images. The technical proof for this
claim is analogous to Step S2 in the proof of Thm. 3 in~\cite{Villanueva2015}.

Finally, observe that the only difference between~\eqref{eq::GDIaux} and
the condition of Corollary~\ref{cor::GGDI1} and is that in the latter, the 
intersection $Y(t)\cap Z(t)$ has been replaced by the tighter set 
$Y(t) \cap \operatorname*{int}(Z(t))$. However, as $Z(t)$ has a nonempty
interior, we have 
\begin{equation*}
\operatorname{cl}(Y(t)\cap\operatorname{int}(Z(t))) = \operatorname{cl}(Y(t)\cap Z(t)) \; ,
\end{equation*}
i.e., the statement of the corollary is not affected if we replace $Z(t)$ in 
the intersection by its interior. This follows from the fact that the supremum of
a continuous function over any bounded set in $\mathbb R^{n_x}$ coincides with the
maximum of the function over the closure of this set.
Thus, we conclude that the statement of the corollary holds.
\hfill\hfill$\diamond$

\subsection{Complete Proof of Theorem~\ref{thm::outerreach}}

The statement of Theorem~\ref{thm::outerreach} can be obtained by a repeated 
application of Corollary~\ref{cor::GGDI1}. First, the corollary is applied to the
reversed differential equation 
 \begin{equation*}
\dot z(t) = -f( z(t),u_1(T-t), u_2(t) ) \; \; 
\text{with} \; \; z(0) \in \mathbb X_2(T) \; .
\end{equation*}
with $Y(t) = Y_{\rm B}(T-t)$, $Z(t) = \mathbb X_2(T-t)$, and $X_0 = \mathbb X_2(T)$.
Reversing time once more, shows that the inequalities for $\dot V[Y_{\rm B}](c)$
and $V[Y_{\rm B}(t)](c)$ in Theorem~\ref{thm::outerreach} imply the inclusion
$X_{\rm B}[u_1](t) \subseteq Y_{\rm B}(t)$, which is valid for all $t\in[0,T]$. 

Now, we apply Corollary~\ref{cor::GGDI1} directly with $Y = Y_{\rm F}$, 
$Z = Y_{\rm B}$ and $X_0 = \{x_0\}$. This yield the inequalities for 
$\dot V[Y_{\rm F}](c)$ and $V[Y_{\rm F}(t)](c)$ in Theorem~\ref{thm::outerreach}, 
implying that $Y_{\rm F}$ is an enclosure for $X_{\rm F}[u_1]$ on $[0,T]$.
Since  $X_{\rm B}[u_1]$ is, by definition, also an enclosure for $X_{\rm F}[u_1]$
on $[0,T]$, it follows that the inclusion 
$Y_{\rm F}(t)\cap Y_{\rm B}(t) \supseteq X_{\rm F}[u_1](t)$ holds for all 
$t\in[0,T]$, yielding the statement of the theorem.
\hfill\hfill$\diamond$

\section{\label{app::ellreach} Proof of Theorem~\ref{thm::ellreach}}

This appendix is divided into a number of subsections which build upon each other
and lead to the proof of Theorem~\ref{thm::ellreach}. Moreover, we 
use the technical convention from Section~\ref{subsec::prelim}.

\subsection{Support functions of set propagation operators}
We recall that the constrained set-propagation operator $\Pi$ has been
introduced in~\eqref{eq::PI}, see Appendix~\ref{app::igdi}, while the shorthand
$\Gamma$ is defined in~\eqref{eq:gammaDef}.

\begin{lem}
\label{lem::dGDI}
Let Assumptions~\ref{ass::rhs1} and~\ref{ass::admissible} be satisfied.
Let $Y,Z:~[0,T] \to \mathbb K_C^{n_x}$ be any set-valued 
functions, such that $Y(t)$ and $Y(t)\cap Z(t)$ are, for all $t\in[0,T]$, strictly
convex and the intersection is nonempty; and such that $V[Y(\cdot)](c)$ 
and $V[Y(\cdot)\cap Z(\cdot)](c)$ are, for all $c \in \mathbb R^{n_x}$,
differentiable. If there exists a continuous function 
$\alpha: \mathbb R \to \mathbb R$ with $\alpha(0)=0$ such that the inequality
\begin{equation*}
V[\Pi(t,t+h,Y(t),Z)](c) \leq V[Y(t+h)](c) + h \alpha(h) \; ,
\end{equation*}
holds for all $c \in \mathbb R^{n_x}$ with $\Vert c \Vert \leq 1$ and all 
$t \in [0,T]$, then
\begin{equation*}
\dot V[Y(t)](c) \geq V[ \Gamma( u_1(t), c, Y(t), Z(t) )](c)
\end{equation*}
holds for all $c \in \mathbb R^{n_x}$ with $\Vert c \Vert \leq 1$ and all 
$t\in[0,T]$.
\end{lem}

\textbf{Proof.}
We first show that the set-propagation operator $\Pi$ satisfies, for all 
$c \in \mathbb R^{n_x}$, the differential inequality
\begin{equation}
\label{eq::auxIneq1}
\begin{aligned}
&\dot V[\Pi(t,t+h,Y(t),Z)](c) \\
&\qquad \qquad \qquad \geq V[ \Gamma( u_1(t), c, Y(t), Z(t) )](c) \; .
\end{aligned}
\end{equation}
The proof of this statement is indirect. 
Let $c \in \mathbb R^{n_x}$ be a vector for which~\eqref{eq::auxIneq1} does
not hold. Then, there exists a point
\begin{equation*}
\xi^{\star} \in \operatorname*{argmax}_{x \in Y(t)} \  c^\intercal x \quad\text{with}
\quad \xi^{\star} \in \operatorname{int}(Z(t)) \;.
\end{equation*}
Otherwise we have $V[ \Gamma( u_1(t), c, Y(t), Z(t) )](c) = -\infty$,
and the inequality~\eqref{eq::auxIneq1} holds.
Now, it follows from the definition of $\Gamma( u_1(t), c, Y(t), Z(t) )$ that
there exists a $\nu_2^{\star} \in \mathbb U_2$ such that
\begin{equation*}
V[ \Gamma( u_1(t), c, Y(t), Z(t) )](c) = 
c^\intercal f\left( \xi^{\star}, u_1(t), \nu_2^{\star} \right) \; ,
\end{equation*}
i.e., we have
\begin{equation*}
 \dot V[\Pi(t,t+h,Y(t),Z)](c) < c^\intercal f\left( \xi^{\star}, u_1(t),\nu_2^{\star} \right) \; .
\end{equation*}
Since $\xi^{\star}\in\operatorname{int}(Z(t))$, the inequality contradicts the 
definition of $\Pi$. Thus,~\eqref{eq::auxIneq1} must hold for all 
$c\in\mathbb{R}^{n_x}$. 

This means that there exists a continuous function $\beta: \mathbb R \to \mathbb R$ such that
\begin{equation}
\label{eq::ineqAux1}
\begin{aligned}
& V[\Pi(t,t+h,Y(t),Z)](c) - V[Y(t)](c) \\
&\qquad \qquad  \geq h V[ \Gamma( u_1(t), c, Y(t), Z(t) )](c) - h \beta(h) \;.
\end{aligned}
\end{equation}
Thus, using the assumptions of this lemma, we can conclude that
\begin{equation*}
\begin{aligned}
&V[Y(t+h)](c)  - V[Y(t)](c) \\
&\qquad\geq V[\Pi(h,Y(t),Z)](c) - V[Y(t)](c) - h \alpha(h) \\
&\qquad\geq h V[ \Gamma( u_1(t), c, Y(t), Z(t) )](c) - h [\alpha(h)+\beta(h)] \; .
\end{aligned}
\end{equation*}
Dividing the last inequality by $h$ on both sides and taking the limit for 
$h \to 0$, we obtain the statement of the lemma.\hfill\hfill$\diamond$

\subsection{Ellipsoidal calculus}
The following proposition summarizes two known results from the field of 
ellipsoidal calculus.
\begin{prop}
\label{prop::intersection}
Let $q_1,q_2 \in \mathbb R^{n_x}$ and $Q_1,Q_2 \in \mathbb S_{++}^{n_x}$ be given.
\begin{enumerate}
\item If $\lambda \in (0,1)$, then 
\begin{equation*}
\mathcal{E}(q_1,Q_1) \oplus \mathcal{E}(q_2,Q_2) \subseteq \mathcal{E} \left(q_1+q_2, \frac{Q_1}{\lambda} + \frac{Q_2}{1-\lambda} \right) \;.
\end{equation*}
\item If $\kappa = ( \kappa_1, \kappa_2 ) \in \mathbb R_+^{2}$ satisfies
\begin{gather*} 
1 =  \kappa_{1}(1-q_{1}^\intercal Q_{1}^{-1}q_{1}) 
+ \kappa_{2}(1-q_{2}^\intercal Q_{2}^{-1}q_2) \\
 + q(\kappa)^\intercal Q(\kappa) q(\kappa) \;,
\end{gather*}
with $Q(\kappa)\in\mathbb{S}^{n_x}_{++}$ and $q(\kappa)\in\mathbb{R}^{n_x}$ given by
\begin{align*}
\tilde Q(\kappa) &= \left( \kappa_{1}Q_{1}^{-1} + \kappa_{2}Q_{2}^{-1} \right)^{-1} \\
\tilde q(\kappa) &= Q_{N}(\kappa) \left(  \kappa_{1}Q_{1}^{-1}q_{1} + 
\kappa_{2}Q_{2}^{-1}q_{2} \right)\;,
\end{align*}
then 
$\mathcal{E}(q_1,Q_1) \cap \mathcal{E}(q_2,Q_2) \subseteq 
\mathcal{E}\left( \tilde q(\kappa), \tilde Q(\kappa) \right)$.
\end{enumerate}
\end{prop}

\textbf{Proof.}
The proofs for these two statements can be found in \cite{Kurzhanski1997}. 
See also~\cite{Houska2011b} for alternative derivations.\hfill\hfill$\diamond$

\begin{cor}
\label{cor::MainCor}
Let Assumptions~\ref{ass::rhs1},~\ref{ass::admissible}, and~\ref{ass::closed}
be satisfied and let $r:\mathbb{R}\to\mathbb{R}^{n_x}$
and $R:\mathbb{R}\to\mathbb{S}^{n_x}_{++}$ be given differentiable functions.
Let $q: \mathbb{R} \to \mathbb{R}^{n_x}$ and $Q: \mathbb{R} \to \mathbb{S}_{++}^{n_x}$
be differentiable and satisfy the differential equations
\begin{align*}
\dot q(t) &= f( q(t),u(t),v ) + \varphi_{3}(q(t),r(t),Q(t),R(t),\kappa(t)) \\
\dot Q(t) &= \Phi_{1}(Q(t),A(t))  + \Phi_{2}\left(Q(t),B(t)VB(t)^\intercal,\sigma(t)\right) \notag \\
& \hphantom{{}={}} + \Phi_{2}\left(Q(t),\Omega(A(t),B(t),q(t),u_1(t),v,Q(t)),\mu(t) \right) \notag \\
& \hphantom{{}={}} + \Phi_{3}\left(q(t),r(t),Q(t),R(t),\kappa(t) \right) \; ,
\end{align*}
on the interval $[0,T]$, for any given functions $A: \mathbb{R}\to \mathbb R^{n_x \times n_x}$, 
$B: \mathbb{R} \to \mathbb R^{n_x \times n_u}$, $\sigma,\mu: \mathbb{R} \to \mathbb R_{++}$
and $\kappa: \mathbb{R} \to \mathbb R_{+}$. Then, the ellipsoidal set valued 
function $Y$ with $Y(t) = \mathcal E( q(t), Q(t) )$ satisfies, 
for all $c \in \mathbb R^{n_x}$, all $t \in [0,T]$, and any 
$Z:\mathbb{R}\to\mathbb{K}^{n_x}_{\rm C}$ with $Z(t) \subseteq \mathcal{E}(r(t),R(t))$
the differential inequality
\begin{align*}
\dot V[Y(t)](c) &\geq V[ \Gamma( u_1(t), c, Y(t), Z(t) )](c) \\
V[Y(t)](c) &\geq V[X_0](c)
\end{align*}
with $X_0 = \mathcal E( q(0), Q(0) )$.
\end{cor}

\textbf{Proof.}
This proof relies on the application of Lemma~\ref{lem::dGDI} with 
$Y(t) = \mathcal E(q(t),Q(t))$. Notice that by differentiability of $q$ and $Q$, the
function $V[\mathcal E(q((\cdot),Q(\cdot))](c)$ is also differentiable on $[0,T]$,
for all $c\in\mathbb{R}^{n_x}$. Let
\begin{align*}
\tilde{x}(t,h) =& \hphantom{{}+{}}x(t) + h f( x(t), u_1(t), u_2(t) ) \\
=& \hphantom{{}+{}} x(t) + h f( q(t), u_1(t), v ) + h A(t) (x(t)-q(t)) \\
 & + hB(t)(u_2(t)-v) +h n(t)
\end{align*}
denote an Euler approximation of the original ODE at time $t$, where
$$n(t) \in \mathcal E( 0, \Omega(A(t),B(t),q(t),u_1(t),v(t),Q(t)))$$
and $x(t) \in \mathcal E(q(t),Q(t))$. The second statement of
Proposition~\ref{prop::intersection} implies that by setting
\begin{equation*}
\begin{aligned}
\tilde Q(t,h) &= \left( \kappa_1(t,h) Q(t)^{-1} + \kappa_2(t,h) R(t)^{-1} \right)^{-1} \\
\tilde q(t,h) &= \left( \kappa_1(t,h) Q(t)^{-1} q(t)  + \kappa_2(t,h) R(t)^{-1}r(t) \right) \tilde Q(t,h) \;.
\end{aligned}
\end{equation*}
we have, by our assumption $Z(t) \subseteq \mathcal{E}(r(t),R(t))$, 
\begin{equation*}
Y(t) \cap Z(t) \subseteq \mathcal E( \tilde q(t,h), \tilde Q(t,h) )
\end{equation*}
as long as $\kappa_1(t,h),\kappa_2(t,h) \geq 0$ satisfy
\begin{equation}
\label{eq::Kconstraint}
\begin{aligned}
\tilde q(t,h)^\intercal \tilde Q(t,h) \tilde q(t,h) =& \hphantom{{}-{}} 1 - \kappa_{1}(t,h)(1-q(t)^\intercal Q(t)^{-1}q(t))  \\
 &- \kappa_{2}(t+h)(1-r(t)^\intercal R(t)^{-1}r(t))\;.
\end{aligned}
\end{equation}
Moreover, since $\mathbb U_2 \subseteq \mathcal E( v, V )$, a repeated application of the first statement in Proposition~\ref{prop::intersection} shows that setting
\begin{eqnarray*}
\hat Q(t,h) &=& \frac{1}{\lambda_1(t,h)} (I + hA(t)) \tilde Q(t,h) (I + hA(t))^\intercal \\
& & + \frac{h^2}{\lambda_2(t,h)} B(t) V B(t)^\intercal \\
& & + \frac{h^2}{\lambda_3(t,h)} \Omega(A(t),B(t),q(t),u_1(t),v,Q(t)))  \\
\hat q(t,h) &=& q_{N}(t+h) + h f( q_{N}(t+h), u_1(t), v)\; , 
\end{eqnarray*}
for any function $\lambda_1(t,h), \lambda_2(t,h), \lambda_3(t,h) > 0$ with $\lambda_1(t,h) + \lambda_2(t,h) + \lambda_3(t,h) = 1$, implies
\begin{align}
\label{eq::piBound}
\begin{aligned}
&V[\Pi(t+h,t, Y(t), Z )](c) \\
&\qquad\quad\leq V[\mathcal E( \hat q(t,h), \hat Q(t,h) )](c)  + h \gamma(h)
\end{aligned}
\end{align}
for a continuous function $\gamma: \mathbb R \to \mathbb R$ with $\gamma(0)=0$, 
for all $c$ with $\Vert c \Vert = 1$, because the Euler discretization is accurate for $h \to 0$.
Next, we substitute
\begin{equation*}
\begin{aligned}
\lambda_1(t,h) &= 1 - h \sigma(t) - h \mu(t)  \\
\lambda_2(t,h) &= h \sigma(t)  \\
\lambda_3(t,h) &= h \mu(t)  \\
\kappa_1(t,h) &= 1 - h \kappa(t) \; ,
\end{aligned}
\end{equation*}
while $\kappa_2(t,h)$ is defined implicitly by~\eqref{eq::Kconstraint}. Differentiating the
above formulas for $\hat q$ and $\hat Q$ with respect to $h$ yields
\begin{align*}
\frac{\rm d}{{\rm d}h} \hat q(t,0) = &\hphantom{{}+{}} 
f\left( \hat q(t,0),u_1(t),v \right) \\
&+ \varphi_{3}\left(\hat q(t,0),r(t),\hat Q(t,0),R(t),\kappa(t) \right) \\
\frac{\rm d}{{\rm d} h}\hat Q(t,0) 
= &\hphantom{{}+{}} \Phi_{1}\left(\hat Q(t,h),A(t)\right)  \\
&+ \Phi_{2}\left(\hat Q(t,0),B(t)VB(t)^\intercal,\sigma(t)\right) \\
&+ \Phi_{2}\left(\hat Q(t,0),\Omega(A(t),B(t),
\widetilde q(t,t), \right .\\
&\qquad \qquad \qquad \left. u_1(t),v,\hat Q(t,0)),\mu(t) \right) \\
&+ \Phi_{3}\left(\hat q(t,0),r(t),\hat Q(t,0),R(t),
\kappa(t) \right) \; .
\end{align*}
Since the right-hand of these derivatives coincide with the differential equations for $q$ and $Q$, we must have
\begin{equation*}
q(t+h) = \hat q(t,h) + O(h^2) \quad \text{and} \quad Q(t,h) = \hat Q(t,0) + O(h^2) \; ,
\end{equation*}
i.e.,~\eqref{eq::piBound} implies that there exists a continuous function $\alpha: \mathbb R \to \mathbb R$ with $\alpha(0) = 0$ and 
\begin{align}
\label{eq::piBound2}
\begin{aligned}
&V[\Pi(t+h,t, Y(t), Z )](c) \\
&\qquad\quad\leq V[\mathcal E( q(t+h), Q(t+h) )](c)  + h \alpha(h) \; .
\end{aligned}
\end{align}
Thus, the statement of this corollary turns into an immediate of consequence of Lemma~\ref{lem::dGDI}. \hfill\hfill$\diamond$

\subsection{Complete Proof of Theorem~\ref{thm::ellreach}}

The statement of Theorem~\ref{thm::ellreach} follows by 
applying the result of Corollary~\ref{cor::MainCor} twice. Firstly, we apply the 
corollary to the reverse dynamic system
\begin{equation*}
\frac{\partial}{\partial t} z(t) = -f( z(t),u_1(T-t), u_2(t) ) \; \; \text{with}
 \; \; z(0) \in \mathbb X_{2}(T) \; .
\end{equation*}
with $q = q_{\rm B}$, $Q = Q_{\rm B}$, $r = s$, $R = S$ and $q_{\rm B}(0) = s(T)$
and $Q_{\rm B}(0) = S(T)$. This yields an enclosure for the backward tube 
$X_{\rm B}[u_1]$. And secondly, we apply Corollary~\ref{cor::MainCor} to~\eqref{eq::odes}
with $q = q_{\rm F}$, $Q = Q_{\rm F}$, $r=q_{\rm B}$, $R = Q_{\rm B}$, 
$q_{\rm F}(0) = x_0$, and $Q_{\rm F}(0) = 0$, which yields the enclosure of the 
forward tube $X_{\rm F}[u_1]$. Since both 
$Y_{\rm B}(\cdot) = \mathcal{E}(q_{\rm B}(\cdot),Q_{\rm B}(\cdot))$ and 
$Y_{\rm F}(\cdot) = \mathcal{E}(q_{\rm F}(\cdot),Q_{\rm F}(\cdot))$ are enclosures
of $X[u_1]$ on $[0,T]$, constructing a set-valued function $Y_{{\rm F}\cap{\rm B}}$ 
by taking their pointwise-in-time intersection yields the statement of the theorem.\hfill\hfill$\diamond$

\end{document}